\documentstyle{article}

\begin{document}

\begin{center}
{\sc On Generalized Clifford Algebra $C_4^{(n)}$ and $GL_q(2;C)$ Quantum
Group}

A.K.Kwa\'sniewski

\medskip

{\em Bia{\l}ystok University,\\Institute of Informatics, ul.
Sosnowa 64\\PL - 15-887 BIA{\L}YSTOK, POLAND\\College of Applied
Mathematics and Computer Science,\\Czysta 11 PL-15-463
Bia{\l}ystok, POLAND\\e-mail: KWANDR@noc.uwb.edu.pl }
\end{center}

\begin{abstract}

The non commuting matrix elements of matrices from quantum group
$GL_q(2;C)$ with $q\equiv \omega $ being the $n$-th root of unity
are given a representation as operators in Hilbert space with help
of $C_4^{(n)}$
generalized Clifford algebra generators appropriately tensored with unit $%
2\times 2$ matrix infinitely many times. Specific properties of
such a representation are presented . Relevance of generalized
Pauli algebra to azimuthal quantization of angular momentum al\`a
L\'evy -Leblond [10] and to polar decomposition of $SU_q(2;C)$
quantum algebra ala Chaichian and Ellinas [6] is also commented.

The case of $q\in C$, $\left| q\right| =1$ may be treated
parallely.

\end{abstract}

\subsection*{1. Introduction}

In 1989 - on the occasion of 70-th birthday of Luigi Radicati - the authors
of [1] have written a contribution entitled: {\it ''Properties of Quantum} $%
2\times 2${\it \ Matrices''}. They follow there the approach of Ludwig D.
Faddeev et all. to the so-called quantum groups [2,3].

Our note deals with the quantum group $GL_q(2;C)$ which is defined [1-6] via
imposing quantization relations on the matrix elements of $GL(2;C)$ matrices
and correspondingly on $SL(2;C)$ group [1].

The paper [1] of Vocos, Wess and Zumino is at the same time the transparent
illustration of the fact that basic representation of a quantum groups
posses special properties.

Here we provide a construction of basic representation of the quantum group $%
GL_q(2;C)$ with $q\equiv \omega \equiv {\rm exp} \left\{ \frac{2\pi i}n\right\} $
being the $n$-th root of unity - in terms of $C_4^{(n)}$ generalized
Clifford algebra generators` tensor products (see [7] and references
therein).

\noindent Quantization relations [1-6] imposed on the matrix elements of $%
GL(2;C)$ - result immediately right from the construction.

Specific properties of such a representation are studied following [1].
Generalized Pauli algebras appear naturally in azimuthal quantization of
angular momentum al\`a L\'evy -Leblond [10] and also in polar decomposition
of the corresponding deformed algebra $SU_q(2;C)$ i.e. ''quantum algebra''
al\`a Chaichian and Ellinas [6]. This is also to be commented.

\noindent Generalized Pauli algebras appear naturally also in $Z_n$-Quantum
Mechanics [9], in $q$-deformed Heisenberg algebras [14] , [9] as well as in $%
q=\omega $ - deformed quantum oscillator [15], [9] as suggested by L. C.
Biedenharn in [12].

\subsection*{2. Quantization relations}

Let $\left(
\begin{array}{cc}
a & b \\
c & d
\end{array}
\right) \in GL(2;C)$.
Then quantization relations [1-6] have the form:
$$
\begin{array}{ccccccc}
{\rm \{rows\}} &  & \rightarrow &  & ab=q\,ba\,\,\,; &  & cd=q\,dc\,\,\,; \\
&  &  &  &  &  &  \\
{\rm \{columns\}} &  & \downarrow &  & ac=q\,ca\,\,\,; &  & bd=q\,db\,\,\,.
\end{array}
\eqno{(2.1)}
$$

In addition to (2.1) also ''diagonal'' quantization relations are imposed;
these being motivated [1] by the obvious requirement that the product of two
quantized matrices (i.e. those satisfying (2.1) and perhaps something
else..) should be a matrix of the kind [1-6]. Namely if apart from matrix $%
A=\left(
\begin{array}{cc}
a & b \\
c & d
\end{array}
\right) \in GL_q(2;C)$ we are given $A^{\prime }=\left(
\begin{array}{cc}
a^{\prime } & b^{\prime } \\
c^{\prime } & d^{\prime }
\end{array}
\right) \in GL_q(2;C)$ where $a^{\prime },b^{\prime },c^{\prime },d^{\prime
} $ {\bf commute} with $a,b,c,d$ then we expect the matrix $A^{\prime \prime
}=AA^{\prime }=\left(
\begin{array}{cc}
a^{\prime \prime } & b^{\prime \prime } \\
c^{\prime \prime } & d^{\prime \prime }
\end{array}
\right) $ to be also of the same kind i.e. we expect noncommuting quantities
$a^{\prime \prime },b^{\prime \prime },c^{\prime \prime },d^{\prime \prime }$
to satisfy (2.1) and perhaps something else if necessary. For that to do in
addition to (2.1) one requires [1]:

$$
{\rm \{diagonals\}}{\quad \quad \quad \quad \quad }bc=cb\,;\quad
ad-da=(q-\frac 1q)bc\,.\eqno{(2.2)}
$$

\noindent If this is accepted then the noncommuting matrix entries of

\[
A^{\prime \prime }=AA^{\prime }=\left(
\begin{array}{cc}
a^{\prime \prime } & b^{\prime \prime } \\
c^{\prime \prime } & d^{\prime \prime }
\end{array}
\right) =\left(
\begin{array}{cc}
aa^{\prime }+bc^{\prime } & ab^{\prime }+bd^{\prime } \\
ca^{\prime }+dc^{\prime } & cb^{\prime }+dd^{\prime }
\end{array}
\right)
\]
also do satisfy (2.1) and (2.2) with the value of $q$ unchanged. In this
connection let us recall [1] that $A^k\equiv \left(
\begin{array}{cc}
a_k & b_k \\
c_k & d_k
\end{array}
\right) \in GL_{q^2}\left( 2;C\right) $ i.e. the noncommuting matrix entries
of $A^k$ satisfy (2.1) and (2.2) with the value of $q^k$. (It was proven in
[1] that $k$ might be any real number).

\medskip\

\noindent \underline{Commentary:}\quad \quad Conditions (2.1) mean that we
have four pairs spanning ''$q$-quantum planes'' - four copies of these $q$%
-planes.

\medskip\

\noindent \underline{Question:}\quad \quad Is it not then enough to impose
(2.1) only?

Perhaps (2.2) might be a technical consequence of the $q$-geometrical (2.1)
requirements` representation.

Out of (2.1) and (2.2.) quantization defining commutation relations the
authors of [1] derive many interesting properties of $SL_q(2;C)$ quantum
group.

\medskip\

It is our aim here to give a simple construction of $GL_q(2;C)$ quantum
group and consequently - $SL_q(2;C)$ quantum group for the case of $q\equiv
\omega \equiv {\rm exp} \left\{ \frac{2\pi i}n\right\} $. Special cases of such
quantum algebras were considered earlier, of course; see [6].

\subsection*{3. The representation of $GL_\omega (2;C)$ elements with help $%
C_4^{\left( n\right) }$ generators.}

Let $\left\{ \gamma _i\right\} _1^4$ be the set of generators of $%
C_4^{\left( n\right) }$ algebra [7] i.e.
$$
\gamma _i\gamma _j=\omega \gamma _j\gamma _i\,;\quad i<j\,;\quad \gamma
_i^n=id.\quad i,j=1,2,3,4.\eqno{(3.1)}
$$
Let us define now two pairs of tensor products of these generators ($%
x,y,X,Y\in C$):

$$
\begin{array}{cc}
a\equiv x\sigma _1=x\left( \gamma _1\otimes \gamma _3\right) \otimes
I\otimes I\otimes ...\,\,\,; & b\equiv y\sigma _2=y\left( \gamma _2\otimes
\gamma _3\right) \otimes I\otimes I\otimes ...\,\,\,; \\
c\equiv X\Sigma _1=X\left( \gamma _1\otimes \gamma _4\right) \otimes
I\otimes I\otimes ...\,\,\,; & d\equiv Y\Sigma _2=Y\left( \gamma _2\otimes
\gamma _4\right) \otimes I\otimes I\otimes ...\,\,\,.
\end{array}
\eqno{(3.2)}
$$

These are entries of the matrix $A$ from [1] realized as operators in
Hilbert space
\[
\left(
\begin{array}{cc}
a & b \\
c & d
\end{array}
\right) =A\equiv \Gamma \equiv \left(
\begin{array}{cc}
x\sigma _1 & y\sigma _2 \\
X\Sigma _1 & Y\Sigma _2
\end{array}
\right) \,\,;
\]
while the $\Gamma $ - notation underlines the fact that $\sigma _1,\sigma _2$%
\thinspace \ ;\thinspace \ $\sigma _1,\Sigma _1\quad $and$\quad \sigma
_2,\Sigma _2$\ \thinspace ;\ $\,\Sigma _1,\Sigma _2\;$ pairs are obtained
from generators of the four corresponding isomorphic copies of $C_2^{\left(
n\right) }$ generalized Clifford algebra which is called - in this case - a
generalized Pauli algebra [7]. ({The generalized Pauli algebra [7] }$%
C_2^{\left( n\right) }${\ generators were already used to provide a
representation of the Heisenberg commutation relations for the finite group }%
${Z}_{{n}}${\ in [8] while describing Heisenberg modules for non-commutative
two-Torii (see [9] for further connotations)}.

\quad \quad One now easily verifies that the following commutation relations
hold:

$$
\begin{array}{lllll}
{\rm \{rows\}} &  & \sigma _1\sigma _2=\omega \sigma _2\sigma _1 & ; &
\Sigma _1\Sigma _2=\omega \Sigma _2\Sigma _1 \\
{\rm \{columns\}} &  & \sigma _1\Sigma _1=\omega \Sigma _1\sigma _1 & ; &
\sigma _2\Sigma _2=\omega \Sigma _2\sigma _2
\end{array}
\eqno{(3.3)}
$$

altogether with
$$
\sigma _1^n=\Sigma _1^n=\sigma _2^n=\Sigma _2^n=id.\eqno{(3.4)}
$$

Thus we have four $\sigma _1,\sigma _2\,\;;\,\;\sigma _1,\Sigma _1\quad $and$%
\quad \sigma _2,\Sigma _2\;\,;\;\,\Sigma _1,\Sigma _2\;$ pairs - four $%
\omega $-frames of ''$\omega $-quantum space'' representations.

\noindent This is represented by the following pictogram matrix:
\[
\left(
\begin{array}{ccc}
\sigma _1 & \rightarrow & \sigma _2 \\
\downarrow &  & \downarrow \\
\Sigma _1 & \rightarrow & \Sigma _2
\end{array}
\right) \,\,,
\]
showing the order of these $\omega $-commuting entries of quantum matrix $%
\Gamma $ in (2.1) commutation relations.

Due to (3.3) $a,b,c,d\;\;$ $\omega $-commuting entries defined by (3.2)
satisfy (2.1) $q$-mutation relations ''automatically'' i.e\thinspace
\thinspace $ab=q\,ba$\ ;\ $cd=q\,dc$\ ;\ $ac=q\,ca$\ ;\ $bd=q\,db$\ ;\ $%
q=\omega $.

Also $bc=cb$ relation is satisfied ''automatically'' i.e. all is due to the
representation .

\medskip\

As for the commutation relations - in order to be complete - we have to
consider also the diagonal directions:
\[
\left(
\begin{array}{ccc}
\sigma _1 & \rightarrow & \sigma _2 \\
\downarrow & \searrow\hskip -3mm\swarrow & \downarrow \\
\Sigma _1 & \rightarrow & \Sigma _2
\end{array}
\right) \,\,.
\]

This is an easy task and one readily verifies that
$$
{\rm \{diagonals\}}\quad \quad \quad \sigma _2\bullet \Sigma _1=\Sigma
_1\bullet \sigma _2\,\,;\;\sigma _1\bullet \Sigma _2=\omega ^2\Sigma
_2\bullet \sigma _1\,\,.\eqno{(3.5)}
$$
In this connection note that $\sigma _1\bullet \Sigma _2=\omega ^2\Sigma
_2\bullet \sigma _1$ is equivalent to $\sigma _1\bullet \Sigma _2-\Sigma
_2\bullet \sigma _1=\left( \omega -\overline{\omega }\right) \sigma _2\Sigma
_1=\left( \omega -\frac 1\omega \right) \sigma _2\Sigma _1$. At the same
time the second commutation relation in (3.5) appears in [1] as a demand in
the form:
\[
ad-da=\left( q-\frac 1q\right) bc\,\,,
\]
which for $q=\omega $ might be rewritten as
$$
ad-da=(\omega -\overline{\omega })bc\,.\eqno{(3.6)}
$$
In view of (3.2) the requirement (3.6) imposes the following bound on
co-ordinates of the four $\omega $-quantum planes
$$
xY=yX\eqno{(3.7)}
$$
However because of $\,\,\sigma _1\bullet \Sigma _2=\omega ^2\Sigma _2\bullet
\sigma _1$, which is equivalent to \thinspace \thinspace $\sigma _1\bullet
\Sigma _2-\Sigma _2\bullet \sigma _1=\left( \omega -\overline{\omega }%
\right) \sigma _2\Sigma _1\,\,$ one may be tempted to replace (see (2.2) )
the quantization condition
$$
ad-da=(q-\frac 1q)bc\eqno{(3.8)}
$$
in the case of $q\equiv \omega \equiv {\rm exp} \left\{ \frac{2\pi i}n\right\} $
by the one resulting from representation of $\omega $-frames i.e.
$$
\lbrack a,d]=\left( 1-\omega ^2\right) \quad \equiv \quad ad=\omega ^2da\,.%
\eqno{(3.9)}
$$
If (3.9) is required instead of (3.8) then we do not have restriction (3.7)
on the co-ordinates of the four $\omega $-quantum planes.

\noindent However then, one should check whether the product of two
quantized matrices (i.e. those satisfying (2.1), $bc=cb$ and (3.9) ) gives a
matrix satisfying the same $q$-mutation relations.

Investigation along lines of [1] is plausible and is now being carried out.
In all formulas above one may take $\omega $ ( condition (3.4) being of
course rejected)) to be $\omega ={\rm exp} \left\{ 2\pi i\alpha \right\} $, with $%
\alpha $ irrational. The algebras thus generated by (3.3) are no more
generalized Pauli algebras of the standard type [7] and they are no more
finite dimensional. Nevertheless these algebras deserve to be called
infinite dimensional Pauli algebras. In this case (3.3) also implies (3.5)
and (3.3) is automatically satisfied in the (3.2) representation.

\noindent The representation of $A^{\prime }=\left(
\begin{array}{cc}
a^{\prime } & b^{\prime } \\
c^{\prime } & d^{\prime }
\end{array}
\right) \in GL_q(2;C)$ where $a^{\prime },b^{\prime },c^{\prime },d^{\prime
} $ {\bf commute} with $a,b,c,d$ is the following
$$
\begin{array}{c}
a^{\prime }\equiv x\sigma _1^{\prime }=xI\otimes I\otimes \left( \gamma
_1\otimes \gamma _3\right) \otimes I\otimes I\otimes ...\,\,\,; \\
b^{\prime }\equiv y\sigma _2^{\prime }=yI\otimes I\otimes \left( \gamma
_2\otimes \gamma _3\right) \otimes I\otimes I\otimes ...\,\,\,; \\
c^{\prime }\equiv X\Sigma _1^{\prime }=XI\otimes I\otimes \left( \gamma
_1\otimes \gamma _4\right) \otimes I\otimes I\otimes ...\,\,\,; \\
d^{\prime }\equiv Y\Sigma _2^{\prime }=YI\otimes I\otimes \left( \gamma
_2\otimes \gamma _4\right) \otimes I\otimes I\otimes ...\,\,\,.
\end{array}
\eqno{(3.2)'}
$$
These are entries of the matrix $A^{\prime }$ from [1] realized as operators
in Hilbert space
\[
\left(
\begin{array}{cc}
a^{\prime } & b^{\prime } \\
c^{\prime } & d^{\prime }
\end{array}
\right) =A^{\prime }\equiv \Gamma ^{\prime }\equiv \left(
\begin{array}{cc}
x\sigma _1^{\prime } & y\sigma _2^{\prime } \\
X\Sigma _1^{\prime } & Y\Sigma _2^{\prime }
\end{array}
\right) \,\,;
\]
Thus moving corresponding $\left( \gamma _i\otimes \gamma _j\right) $ by
step two to the right we obtain others $A^{\prime \prime },A^{\prime \prime
\prime }$ etc. apart from their products $AA^{\prime \prime },A^{\prime
\prime }A,...,A^{\prime \prime }A^{\prime \prime \prime }GL_q(2;C)$.

\medskip\

\noindent {\bf General conclusion: }

{\it The bulk of implications} of (2.1) and (2.2) for the quantum group
[1-6] in the case of $q\equiv \omega \equiv {\rm exp} \left\{ \frac{2\pi i}%
n\right\} $ or $\omega ={\rm exp} \left\{ 2\pi i\alpha \right\} $, with $\alpha $
irrational {\it results from the (3.2) representation} of - four $\omega $%
-frames of ''$\omega $-quantum space'' {\it with (3.7) bound on coordinates}
being imposed. The possibility of replacing (3.8) by (3.9) deserves to be
investigated.

\medskip\

Among the bulk of implications of (2.1) and (2.2) let us here only note [1]
that $A^k$ satisfies (2.1) and (2.2) with $q^k$ quantum parameter what
results in our case of $q\equiv \omega $ being the $n$-th root of unity in
cyclic arriving to the $q=1$ case.

\subsection*{4. Interpretation of quantization relations}

Now we shall quote and comment interpretations of quantization relations as
stated by
$$
\begin{array}{ccccccc}
{\rm \{rows\}} &  & \rightarrow &  & ab=q\,ba\,\,\,; &  & cd=q\,dc\,\,\,; \\
&  &  &  &  &  &  \\
{\rm \{columns\}} &  & \downarrow &  & ac=q\,ca\,\,\,; &  & bd=q\,db\,\,\,.
\end{array}
\eqno{(2.1)}
$$
$$
{\rm \{diagonals\}}{\quad \quad \quad \quad \quad }bc=cb\,;\quad
ad-da=(q-\frac 1q)bc\,.\eqno{(2.2)}
$$
We shall follow [1] and then [6] keeping in mind that for $q=\omega $ the
above quantization relations result from the (3.2) representation with (3.7)
bound on coordinates .

\noindent Recall first that the quantum or $q$-deformed determinant of the
matrix $A=\left(
\begin{array}{cc}
a & b \\
c & d
\end{array}
\right) \in GL_q(2;C)$ is defined by
$$
D_q={\rm det}_qA=ad-qbc=da-q^{-1}bc\,\,.\eqno{(4.1)}
$$
$D_q$ commutes with all elements satisfying quantization relations (2.1) and
(2.2). It might be also shown [1] that
$$
{\rm det}_{q^2}A^k=\left( {\rm det}_qA\right) ^k\eqno{(4.2)}
$$
In our case of $q=\omega $ this means that ${\rm det\,}A^n=\left( {\rm det\,}%
_\omega A\right) ^n$; hence $\omega $-deformed determinant of quantum matrix
$A$ is obtainable form usual determinant of a matrix $B=A^n$ with commuting
entries; {\{note that }${A}^{{k}}${\ satisfies (2.1) and (2.2) with }${q}^{{k%
}}\}$.

Let us also [1] introduce the $q$-deformed or quantum epsilon matrix
$$
\varepsilon _q=\left(
\begin{array}{cc}
0 & \frac 1{\sqrt{q}} \\
-\sqrt{q} & 0
\end{array}
\right) \eqno{(4.3)}
$$
Of course $\varepsilon _q^2=-1$ and as one may check it [1] the quantization
relations (2.1) and (2.2) \{{resulting for }${q=\omega }${\ from the (3.2)
representation with (3.7) bound on coordinates\}}are equivalent to
$$
A^T\varepsilon _qA=A\varepsilon _qA^T=D_q\varepsilon _q\eqno{(4.4)}
$$

Therefore for $SL_q(2;C)$ quantum group the \underline{{\bf quantization
relations}} may be interpreted as $q$-deformed or quantum \underline{{\bf %
symplectic conditions}} imposed on matrices $A$ under quantization. Let us
now come over to the other characterization of quantization relations (2.1)
and (2.2) resulting for $q=\omega $ from the (3.2) representation with (3.7)
bound on coordinates.

Namely let us interpret matrices $A=\left(
\begin{array}{cc}
a & b \\
c & d
\end{array}
\right) \in GL_q(2;C)$ as endomorphism acting on a quantum plane with points
labeled by $\left(
\begin{array}{c}
x \\
y
\end{array}
\right) ,\left(
\begin{array}{c}
x^{\prime } \\
y^{\prime }
\end{array}
\right) ,...$ etc where $xy=\omega yx\,,\,\,x^{\prime }y^{\prime }=\omega
y^{\prime }x^{\prime },..$ etc.

Now if $A$ is supposed to map $\omega $-quantum plane onto the same $\omega $%
-quantum plane
$$
\left(
\begin{array}{cc}
a & b \\
c & d
\end{array}
\right) \left(
\begin{array}{c}
x \\
y
\end{array}
\right) =\left(
\begin{array}{c}
x^{\prime } \\
y^{\prime }
\end{array}
\right) \eqno{(4.5)}
$$
then quantization relations( 2.1) and (2.2) are necessary and sufficient
condition for that to be the case [6] , [11].

Coordinates $x$ \& $y$ of $\omega $-quantum plane commute with entries of $A$
matrix as it is the case with product of matrix $A=\left(
\begin{array}{cc}
a & b \\
c & d
\end{array}
\right) \in GL_q(2;C)$ with $A^{\prime }=\left(
\begin{array}{cc}
a^{\prime } & b^{\prime } \\
c^{\prime } & d^{\prime }
\end{array}
\right) \in GL_q(2;C)$ where $a^{\prime },b^{\prime },c^{\prime },d^{\prime
} $ {\bf commute} with $a,b,c,d$.

If entries of $A$ matrix are represented as in (3.2) then there exist
infinitely many representations of $\omega $-quantum plane noncommuting
coordinates; for example
$$
\begin{array}{ll}
x=xI\otimes I\otimes \left( \gamma _1\otimes \gamma _3\right) \otimes
I\otimes I\otimes ...\,\,\,; & y=yI\otimes I\otimes \left( \gamma _2\otimes
\gamma _3\right) \otimes I\otimes I\otimes ...\,\,\,; \\
{\rm or} &  \\
x=xI\otimes I\otimes \left( \gamma _1\otimes \gamma _4\right) \otimes
I\otimes I\otimes ...\,\,\,; & y=yI\otimes I\otimes \left( \gamma _2\otimes
\gamma _4\right) \otimes I\otimes I\otimes ...\,\,\,.
\end{array}
\eqno{(4.6)}
$$
or others obtained by moving corresponding $\left( \gamma _i\otimes \gamma
_j\right) $ to the right.

\medskip\

The observation concerning (4.5) was made in [6] on the occasion of quantum
polar decomposition of the algebra of $SU_q(2)$ quantum group. The polar
decomposition of the $SU(2)$ group algebra was provided by L\'evy-Leblond in
his Mexican paper [10]. There he had interpreted such a polar decomposition
as a tool for\quad ''azimuthal quantization of angular momentum''.

\noindent It appears that in both cases generalized Pauli matrices ({%
building blocks for all representations of generalized Clifford algebras)
appear in the same way [6]. }

\subsection*{5. Polar decomposition of $SU_q(2;C)$ and $SU(2;C)$ groups`
algebras}

The standard basis of Lie algebra $su(2)=so(3)$ is well known to be
represented by:

$$
\begin{array}{c}
J_3=\sum\limits_{m=-j}^jm\left| jm\right\rangle \,\,\,; \\
J_{+}=\sum\limits_{m=-j}^j\sqrt{\left( j-m\right) \left( j+m+1\right) }%
\left| j\left( m+1\right) \right\rangle \,\left\langle jm\right| \,\,\,; \\
J_{-}=\sum\limits_{m=-j}^j\sqrt{\left( j+m\right) \left( j-m+1\right) }%
\left| j\left( m+1\right) \right\rangle \,\left\langle jm\right| \,\,\,.
\end{array}
\eqno{(5.1)}
$$

In [12] Biedenharn proposed a new realization of quantum group $SU_q\left(
2\right) $ and in order to realize generators of a $q$-deformation $%
U_q(su(2))$ of the universal enveloping algebra of the Lie algebra $su(2)$
he defined a pair of mutual commuting $q$-harmonic oscillator systems (al\`a
Jordan-Schwinger approach to $su(2)$ Lie algebra).

At the same time in [13] Mac Farlane had also provided us with $q$%
-oscillator description of $SU_q\left( 2\right) $ (al\`a Jordan-Schwinger
approach to $su(2)$ Lie algebra).

\noindent The generators of a $q$-deformation $U_q(su(2))$ of the universal
enveloping algebra of the Lie algebra $su(2)$ (called by physicists
''generators of the quantum group $SU_q\left( 2\right) $'' - which is not a
group!) are obtained from (5.1) by one of several possible $q$-deformations.
In Biedenharn's and Mac Farlane's case one uses the following $q$%
-deformation of numbers and operators:

$$
\left[ x\right] _q=\frac{q^x-q^{-x}}{q-q^{-1}}\,\,\,.\eqno{(5.2)}
$$
Thus $q$-deformed (5.1) representation of generators now reads
$$
\begin{array}{c}
J_3=\sum\limits_{m=-j}^jm\left| j,m\right\rangle _q\,\,\,; \\
J_{+}=\sum\limits_{m=-j}^j\sqrt{\left[ j-m\right] _q\left[ j+m+1\right] _q}%
\left| j,\left( m+1\right) \right\rangle _{q\,\,q}\left\langle j,m\right|
\,\,\,; \\
J_{-}=\sum\limits_{m=-j}^j\sqrt{\left[ j+m\right] _q\left[ j-m+1\right] _q}%
\left| j,\left( m+1\right) \right\rangle _{q\,\,q}\left\langle j,m\right|
\,\,\,,
\end{array}
\eqno{(5.3)}
$$
where
$$
\left| j,m\right\rangle _q=\left| j+m\right\rangle _q\left| j-m\right\rangle
_q=\frac{{a_{1q}^{+}}^{j+m}{a_{2q}^{+}}^{j-m}}{\left[ j+m\right] _q!\left[
j-m\right] _q!}\left| 0\right\rangle _q\eqno{(5.4)}
$$
and ${a_{1q}^{+}\,}$, ${a_{2q}^{+}}$ represent two mutually commuting
creation operators of $q$-quantum harmonic oscillators. Corresponding
commutation relations of the generators of a $q$-deformation $U_q(su(2))$ of
the universal enveloping algebra of the Lie algebra $su(2)$ are of the
familiar though now $q$-deformed form [12], [13], [6]:
$$
\lbrack J_3,J_{+}]=J_{+}\,\,\,;\quad [J_3,J_{-}]=-J_{-}\,\,\,;\quad
[J_{+},J_{-}]=[2J_3]_q\,\,\,.\eqno{(5.5)}
$$
From (5.3) one may derive [6] the polar decomposition of the generators
$J_{+},J_{-}$:

$$
J_{+}=\sqrt{J_{+}J_{-}}\sigma _1^{-1}=\sigma _1^{-1}\sqrt{J_{-}J_{+}}\quad
\quad J_{-}=\sqrt{J_{+}J_{-}}\sigma _1=\sigma _1\sqrt{J_{-}J_{+}}\,\,\,,%
\eqno{(5.6)}
$$
where $\sigma _1$ is the first of the two generators of generalized Pauli
algebra [7], [9]

$$
\sigma _1=\left(
\begin{array}{ccccccc}
0 & 1 & 0 & 0 & \cdots & 0 & 0 \\
0 & 0 & 1 & 0 & \cdots & 0 & 0 \\
0 & 0 & 0 & 1 &  & 0 & 0 \\
\vdots &  & \vdots &  & \vdots &  & \vdots \\
0 & 0 & 0 & 0 &  & 1 & 0 \\
0 & 0 & 0 & 0 & \cdots & 0 & 1 \\
1 & 0 & 0 & 0 & \cdots & 0 & 0
\end{array}
\right) =\left( n\times n\right) \eqno{(5.7)}
$$

The second generator $\sigma _2$ has been also used in [6] in order to
remark on relevance of such a pair $\sigma _1$ , $\sigma _2$ to $GL_\omega
(2;C)$ properties. It is to be noted here that the polar decomposition for
undeformed $su(2)=so(3)$ algebra of undeformed quantum angular momentum had
been performed already by L\'evy -Leblond in his Mexican paper [10]. There
he had interpreted such a polar decomposition as the ''azimuthal
quantization of angular momentum''. And it should be also noted here -
following the authors of [6] that - what we know as - generalized Pauli
algebra appears in $q$-deformed and in undeformed case of polar
decomposition in the same way (5.6).

Neither L\'evy -Leblond nor the authors of [6] had realized that they are
dealing with generalized Pauli algebra [9]. These has been realized
afterwards by T. S. Santhanam [15].

In the notation of [9] and earlier papers quoted there
$$
\sigma _2\equiv U=\omega ^Q={\rm exp} \left\{ \frac{2\pi i}nQ\right\} =\left(
\begin{array}{cccc}
1 & 0 & \cdots & 0 \\
0 & \omega &  & 0 \\
\vdots &  & \ddots & \vdots \\
0 & 0 & \cdots & \omega ^{n-1}
\end{array}
\right) \eqno{(5.8)}
$$
where
$$
Q=\left(
\begin{array}{cccc}
0 & 0 & \cdots & 0 \\
0 & 1 &  & 0 \\
\vdots &  & \ddots & \vdots \\
0 & 0 & \cdots & n-1
\end{array}
\right) \,\,;\eqno{(5.9)}
$$
$$
\sigma _1\equiv V=\omega ^P={\rm exp} \left\{ \frac{2\pi i}nP\right\} =\left(
\begin{array}{ccccc}
0 & 1 & 0 & \cdots & 0 \\
0 & 0 & 1 & \ddots & \vdots \\
\vdots &  & \ddots & \ddots & 0 \\
0 &  &  & \ddots & 1 \\
1 & 0 & \cdots & 0 & 0
\end{array}
\right) \eqno{(5.10)}
$$
where
$$
P=S^{+}QS=\left( P_{\alpha ,\kappa }\right) \eqno{(5.11)}
$$
\[
P_{\alpha ,\kappa }=\left\{
\begin{array}{ll}
0 & \quad \alpha =\kappa \\
\left[ \overline{\omega }^{\left( \alpha -\kappa \right) }-1\right] ^{-1} &
\quad \alpha \neq \kappa
\end{array}
\right.
\]
and\quad
$$
S=\left( <\widetilde{k}|l>\right) =\frac 1{\sqrt{n}}\left( \omega
^{kl}\right) _{k,l\in Z_n}\eqno{(5.12)}
$$
$\quad $ is the Sylvester matrix.

Formulas (5.8) - (5.12) contain the main information on quantum kinematics
of the finite dimensional quantum mechanics. Here we interpret polar
decomposition of quantum angular momentum algebra $su(2)=so(3)$ formalism as
a model of finite dimensional quantum mechanics with the classical phase
space being the torus $Z_n\times Z_n$ (see [9] and references therein) .
This possibility was already considered by Weyl in [16].

''Azimuthal quantization of angular momentum'' was interpreted afterwards as
the finite dimensional quantum mechanics by Santhanam et. all [15].

The considerations of this section allow us to hope to elaborate soon more
on the $q$-deformed finite dimensional quantum mechanics treated as an
interpretation of $q$-deformed $su(2)$ algebra of $q$-deformed angular
momentum (for example by Schwinger method).

\bigskip\

\noindent {\bf Acknowledgments}

\smallskip\

Discussions with W. Bajguz are highly acknowledged. The author expresses his
gratitude to him also for LaTeX version of my work.

\bigskip\

\end{document}